\documentclass[12pt,twoside]{article}
\usepackage{amsgen,amsmath,amstext,amsbsy,amsopn,amsfonts,amssymb}
\usepackage{jltmac2e}
\renewcommand{\version}{{\rm Version of~}\date}

\newcommand\R{{\mathbb{R}}}

\def\Im{{\rm Im}\,}
\def\Re{{\rm Re}\,}
\def\const{{\rm const.}\,}

\renewcommand{\vol}{??}
\renewcommand{\annum}{??}
\renewcommand{\title}{{
 \centerline{Equipartition of energy for the wave equation}
 \centerline{associated to the Dunkl-Cherednik Laplacian}}}
\renewcommand{\author}{Fatma Ayadi\thanks{\; This work is part of my master thesis, done under the supervision of M.~Sifi (Universit\'e de Tunis El Manar) and J.--Ph.~Anker (Universit\'e d'Orl\'eans), within the program CMCU 07G 1501
 \it Aspects analytique et probabiliste de la th\'eorie de Dunkl}}
\renewcommand{\lastname}{Ayadi}

\def\rec{}
\def\fin{}

\renewcommand{\address}{F. Ayadi

Universit\'e de Tunis El Manar

Facult\'e des Sciences de Tunis

D\'epartement de Math\'ematiques

2092 Tunis El Manar, Tunisie

                        \catcode`\@=11
                        ayfatm@yahoo.fr
                        \catcode`\@=13

\&

Universit\'e d'Orl\'eans \& CNRS

F\'ed\'eration Denis Poisson (FR 2964)

Laboratoire MAPMO (UMR 6628)

B\^atiment de Math\'ematiques

B.P. 6759, 45067 Orl\'eans cedex 2, France

\catcode`\@=11
                        Fatma.Ayadi@etu.univ-orleans.fr
                        \catcode`\@=13}

\begin{document}
\firstpage

\newtheorem{odstavec}[Theorem]{\hskip -2mm}

\begin{abstract}This paper is concerned with energy properties of the wave equation
associated to the Dunkl-Cherednik Laplacian.
We establish the conservation of the total energy,
the strict equipartition of energy under suitable assumptions
and the asymptotic equipartition in the general case.\\
{\sl Mathematics Subject Index 2000:} {Primary 35L05\,; Secondary 22E30, 33C67, 35L65, 58J45}
\end{abstract}

\vbadness=100000

 \section{Introduction}

\noindent
We use \cite{O2} as a reference for the Dunkl-Cherednik theory.
Let $\mathfrak{a}$ be a Euclidean vector  space of dimension $d$ equipped
with an inner product $\langle \cdot,\cdot\rangle.$ Let
$\mathcal{R}$ be a crystallographic root system in $\mathfrak{a},$
$\mathcal{R}^+$ a positive subsystem and
$W$ the Weyl group generated by the reflections
\,$r_\alpha(x)\!=\!x\!
-\!2\frac{\langle\alpha,x\rangle}{\|\alpha\|^2}\alpha$
\,along the roots $\alpha\!\in\!\mathcal{R}.$ We let  $k:\mathcal{R}\to[\,0,+\infty)$
denote a   multiplicity function on the root system $\mathcal R,$ and
$\rho=\frac12\sum_{\alpha\in\mathcal{R}^+}k_\alpha\alpha$.
We note that $k$ is $W$-invariant. The Dunkl-Cherednik operators are
the following differential-difference operators,
which are deformations of partial derivatives
and still commute pairwise\,:
\begin{equation*}
T_\xi f(x)
=\partial_{\xi}f(x)-\langle\rho,\xi\rangle f(x)
+\sum_{\alpha\in\mathcal{R}^{+}}k_\alpha\,
\frac{\langle\alpha,\xi\rangle}{1\,-\,e^{-\langle\alpha,x\rangle}}\,
\{f(x)\!-\!f(r_{\alpha}x)\}\,.
\end{equation*}

Given an orthonormal basis $\{\xi_1,\dots,\xi_d\}$ of $\mathfrak{a}$,
the Dunkl-Cherednik Laplacian is defined by
\begin{equation*}
Lf(x)=\sum_{j=1}^{d}T_{\xi_j}^{2}f(x)\,.
\end{equation*}
More explicit formulas for $L$ exist but they will not be used in
this paper. The Laplacian  $L$ is selfadjoint with respect to the
measure $\mu(x) dx$ where
\begin{equation*}
\mu(x)=\prod_{\alpha\in\mathcal{R}^+}
\bigl|\,2\sinh\frac{\langle\alpha,x\rangle}2\bigr|^{2k_\alpha}\,.
\end{equation*}

Consider the wave equation
\begin{equation}
\label{WaveEquation}
\begin{cases}
\,\partial_t^2u(t,x)=L_xu(t,x),\\
\,u(0,x)=f(x),
\,\partial_t\big|_{t=0}u(t,x)=g(x),
\end{cases}
\end{equation}
with smooth and compactly supported initial data $(f,g).$
Let us introduce:
\begin{itemize}
\item
the \textit{kinetic energy}
\,$\mathcal{K}[u](t)
=\frac12{\displaystyle\int_{\mathfrak{a}}}\,
\big|\partial_tu(x,t)\big|^{2}\mu(x)\,dx$,
\item
the \textit{potential energy}
\,$\mathcal{P}[u](t)
=-\,\frac12{\displaystyle\int_{\mathfrak{a}}}\,
Lu(x,t)\,\overline{u(x,t)}\,\mu(x)\,dx$,
\item
the \textit{total energy}
\,$\mathcal{E}[u](t)
=\mathcal{K}[u](t)
+\mathcal{P}[u](t)$.
\end{itemize}

In this paper we prove
\begin{itemize}
\item
the conservation of the total energy:
\begin{equation}
\label{Conservation}
\textstyle
\mathcal{E}[u](t)=\text{constant},
\end{equation}
\item
the strict equipartition of energy,
under the assumptions that the dimension $d$ is odd
and that all the multiplicities $k_\alpha$ are integers:
\begin{equation}
\label{StrictEquipartition}
\textstyle
\mathcal{K}[u](t)=\mathcal{P}[u](t)=\frac12\,\mathcal{E}[u]
\text{ \,for  }|t|\text{  large}\,,
\end{equation}
\item
 the asymptotic equipartition of energy, for arbitrary $d$ and  $\R^+$-valued multiplicity function $k$:
\begin{equation}
\label{AsymptoticEquipartition}
\textstyle
\mathcal{K}[u](t)\to\frac12\,\mathcal{E}[u]
\text{ \,and \,}
\mathcal{P}[u](t)\to\frac12\,\mathcal{E}[u]
\text{ \,as \,}
|t| \text{ goes to } \infty.
\end{equation}
\end{itemize}

The proofs follow \cite{BOS}
and use the Fourier transform in the Dunkl-Cherednik setting,
which we will recall  in next section. We mention that during the past twenty years,
several works were devoted to
Huygens' principle and equipartition of energy
for wave equations
on symmetric spaces and related settings.
See for instance \cite{BO}, \cite [ch.\,V, \S\;5]{H}, \cite{BOS},
\cite{KY}, \cite{SO}, \cite{BOP1}, \cite{BOP2}, \cite{S}.

\section{Generalized hypergeometric functions
and Dunkl-Cherednik transform}

\noindent
Opdam \cite{O1} introduced the following special functions,
which are deformations of exponential functions
\,$e^{\langle\lambda,x\rangle}$,
\,and the associated Fourier transform.

\begin{theorem}
There exist a neighborhood \,$U$ of \,$0$ in $\mathfrak{a}$
and a unique holomorphic function \,$(\lambda,x)\mapsto G_\lambda(x)$
\,on \,$a_{\mathbb{C}}\!\times\!(\mathfrak{a}\!+\!iU)$ such that
\begin{equation*}
\begin{cases}
\;T_\xi G_\lambda(x)=\langle\lambda,\xi\rangle\,G_\lambda(x)
\quad\forall\;\xi\!\in\!\mathfrak{a}\,,\\
\;G_\lambda(0)=1\,.
\end{cases}
\end{equation*}
Moreover, the following estimate holds on
\,$a_{\mathbb{C}}\!\times\!\mathfrak{a}$\,:
\begin{equation*}
|G_\lambda(x)|\le|W|^{\frac12}\,e^{\|\Re\lambda\|\|x\|}\,.
\end{equation*}
\end{theorem}

\begin{Definition} Let $f$ be a nice function   on
$\mathfrak{a}$, say  $f$ belongs to the space
$\mathcal{C}_c^\infty(\mathfrak{a}) $ of smooth functions on $\mathfrak{a}$
with compact support. Its Dunkl-Cherednik transform  is defined by
\begin{equation*}
\mathcal{F}\!f(\lambda)
=\int_{\mathfrak{a}}f(x)\,G_{-iw_0\lambda}(w_0x)\,\mu(x)\,dx.
\end{equation*}
 Here $w_0$ denotes the longest element in the Weyl group $W.$
\end{Definition}

The involvement of $w_0$ in the definition of $\mathcal F$ is related to
the below skew-adjointness property of the Dunkl-Cherednik operators
with respect to the inner product
\begin{equation*}
\langle f,g\rangle=\int_{\mathfrak{a}}f(x)\,\overline{g(x)}\,\mu(x)\,dx,
\qquad f,g\in\mathcal{C}_c^\infty(\mathfrak{a}).
\end{equation*}

\begin{lemma}
The adjoint of \,$T_\xi$ is $-w_0T_{w_0\xi}w_0$\,:
\begin{equation*}
\langle T_\xi f,g\rangle=\langle f,-w_0T_{w_0\xi}w_0g\rangle\,.
\end{equation*}
\end{lemma}

As an immediate consequence, we obtain:
\begin{corollary}
For every \,$\xi,\lambda\!\in\!\mathfrak{a}$
and $f\!\in\!\mathcal{C}_c^\infty(\mathfrak{a})$,
we have
\begin{equation*}
\mathcal{F}(T_\xi f)(\lambda)
=i\langle\lambda,\xi\rangle\mathcal{F}\!f(\lambda),
\end{equation*}
and therefore
\begin{equation*}
\mathcal{F}(Lf)(\lambda)=-\|\lambda\|^{2}\,\mathcal{F}\!f(\lambda)\,.
\end{equation*}
\end{corollary}

Next we will recall from \cite{O1}   three main results about the Dunkl-Cherednik transform
(see also \cite{O2}). For $R>0,$ let  $\mathcal{C}_R^\infty(\mathfrak{a})$ be
the space of smooth functions on $\mathfrak{a}$
\,vanishing outside the ball
\,$B_R\!=\!\{x\!\in\!\mathfrak{a}\,|\|x\|\!\le\!R\}.$
We let $\mathcal{H}_{R}(\mathfrak{a}_{\mathbb{C}})$ denote  the space of holomorphic
functions  $h$ on the complexification $\mathfrak{a}_{\mathbb{C}}$ of $\mathfrak{a}$
such that, for every integer $N\!>\!0$,
\begin{equation*}
\textstyle
\sup_{\lambda\in\mathfrak{a}_{\mathbb{C}}}
(1\!+\!\|\lambda\|)^N\,
e^{-R\,\|\Im\lambda\|}\,|h(\lambda)|
<+\infty\,.
\end{equation*}
\begin{theorem}
\label{PaleyWiener} {\rm\bf (Paley-Wiener)}
 The transformation
$\mathcal{F}$ is an isomorphism of
$\mathcal{C}_R^\infty(\mathfrak{a})$ onto
$\mathcal{H}_{R}(\mathfrak{a}_{\mathbb{C}}),$ for every $R\!>\!0$.
\end{theorem}

The Plancherel formula and the inversion formula of $\mathcal F$
involve the complex measure \,$\nu(\lambda)\,d\lambda$
\,with density
\begin{equation*}
\textstyle
\nu(\lambda)\,
=\hspace{-1mm}\prod\limits_{\alpha\in\mathcal{R}_0^+}\hspace{-1mm}
\frac
{\Gamma\bigl(i\langle\lambda,\check\alpha\rangle\,
+\,k_\alpha\bigr)\vphantom{\big|}}
{\Gamma\bigl(i\langle\lambda,\check\alpha\rangle\bigr)
\vphantom{\big|}}\,
\frac
{\Gamma\bigl(\frac{i \langle\lambda,\check\alpha\rangle\,
+\,k_{\alpha}}2+\,k_{2\alpha}\bigr)\vphantom{\big|}}
{\Gamma\bigl(\frac{i\langle\lambda,\check\alpha\rangle\,
+\,k_\alpha}2\bigr)\vphantom{\big|}}\,
\frac
{\Gamma\bigl(-\,i\langle\lambda,\check\alpha\rangle\,
+\,k_\alpha\bigr)\vphantom{\big|}}
{\Gamma\bigl(-\,i\langle\lambda,\check\alpha\rangle\,+\,1\bigr)
\vphantom{\big|}}\,
\frac
{\Gamma\bigl(\frac{-i \langle\lambda,\check\alpha\rangle\,
+\,k_{\alpha}}2+\,k_{2\alpha}+\,1\bigr)\vphantom{\big|}}
{\Gamma\bigl(\frac{-i\langle\lambda,\check\alpha\rangle\,
+\,k_\alpha}2\bigr)\vphantom{\big|}}\,,
\end{equation*}
where
\,$\mathcal{R}_0^+\!
=\hspace{-.25mm}\{\alpha\!\in\mathcal{R}^+
|\,\frac\alpha2\!\notin\!\mathcal{R}\,\}$
\,is the set of positive indivisible roots,
\,$\check\alpha\!=\!2\|\alpha\|^{-2}\alpha$
\,the coroot corresponding to $\alpha$,
and \,$k_{2\alpha}\!=\!0$ \,if \,$2\alpha\!\notin\!\mathcal{R}$.
Notice that $\nu$ is an analytic function on $\mathfrak{a},$
with a polynomial growth
and which extends meromorphically to~$\mathfrak{a}_{\mathbb{C}}$.
It is actually a polynomial
if the multiplicity function $k$ is integer-valued
and it has poles otherwise.

\begin{theorem}
{\rm\bf (Inversion formula)}
 There is a constant \,$c_0\!>\!0$ such
that, for every $f\!\in\!\mathcal{C}_c^\infty(\mathfrak{a})$,
\begin{equation*}
f(x)\,=\,c_0\int_{\mathfrak{a}}\,
\mathcal{F}\!f(\lambda)\,G_{i\lambda}(x)\,\nu(\lambda)\,d\lambda\,.
\end{equation*}
\end{theorem}

\begin{theorem}\label{Plancherel}
{\rm\bf (Plancherel formula)}
 For every
$f,g\!\in\!\mathcal{C}_c^\infty(\mathfrak{a})$,
\begin{equation*}
\int_{\mathfrak{a}}\,f(x)\,\overline{g(x)}\,\mu(x)\,dx\,
=\,c_0\int_{\mathfrak{a}}\,\mathcal{F}\!f(\lambda)\,
\widetilde{\mathcal{F}}\hspace{-.25mm}g(\lambda)\,
\nu(\lambda)\,d\lambda\,,
\end{equation*}
where \,$\widetilde{\mathcal{F}}\hspace{-.25mm}g(\lambda)
:=\overline{\mathcal{F}(w_0g)(w_0\lambda)}
=\!{\displaystyle\int_{\mathfrak{a}}}\;
\overline{g(x)}\,G_{i\lambda}(x)\,\mu(x)\,dx$\,.
\end{theorem}

\section{Conservation of energy}
\noindent
This section is devoted to the proof of \eqref{Conservation}.
Via the Dunkl-Cherednik transform,
the wave equation \eqref{WaveEquation} becomes
\begin{equation*}
\begin{cases}
\,\partial_t^2\mathcal{F}\hspace{-.25mm}u(t,\lambda)
=-\|\lambda\|^2\mathcal{F}\hspace{-.25mm}u(t,\lambda),\\
\,\mathcal{F}\hspace{-.25mm}u(0,\lambda)
=\mathcal{F}\!f(\lambda),
\,\partial_t\big|_{t=0}\mathcal{F}\hspace{-.25mm}u(t,\lambda)
=\mathcal{F}\!g(\lambda),
\end{cases}
\end{equation*}
and its solution satisfies
\begin{equation}\label{SolutionFourier}
\textstyle
\mathcal{F}\hspace{-.25mm}u(t,\lambda)
=(\cos t\|\lambda\|)\,\mathcal{F}\!f(\lambda)
+\frac{\sin t\|\lambda\|}{\|\lambda\|}\mathcal{F}\!g(\lambda)\,.
\end{equation}
By means of the Paley-Wiener Theorem 2.5, in \cite[p. 52-53]{S} the author proves
 the following finite speed propagation property:
\begin{quote}
Assume that the initial data $f$ and $g$ belong to
$\mathcal{C}_R^\infty(\mathfrak{a})$.
Then the solution $u(t,x)$ belongs to
$\mathcal{C}_{R+|t|}^{\infty}(\mathfrak{a})$
as a function of $x$.
\end{quote}

Let us express the potential and kinetic energies defined in the introduction
via the Dunkl-Cherednik transform. Using the Plancherel formula and Corollary 2.4,
we have
\begin{equation}
\label{PotentialFourier1}
\mathcal{P}[u](t)
=\,{\textstyle\frac{c_0}2}{\displaystyle\int_{\mathfrak{a}}}\;
\|\lambda\|^2\mathcal{F}\hspace{-.25mm}u(t,\lambda)\,
\widetilde{\mathcal{F}}\hspace{-.25mm}u(t,\lambda)\,
\nu(\lambda)\,d\lambda\,.
\end{equation}
Moreover, since the Dunkl-Cherednik Laplacian is $W$-invariant, it follows that  $(w_0u)(t,x)\!=\!u(t,w_0x)$ is the solution to the wave equation \eqref{WaveEquation}
with the initial data \,$w_0f$ and $w_0g$.
Thus
\begin{equation}
\label{SolutionFourierBis}
\textstyle
\widetilde{\mathcal{F}}\hspace{-.25mm}u(t,\lambda)
=(\cos t\|\lambda\|)\,
\widetilde{\mathcal{F}}\!f(\lambda)
+\frac{\sin t\|\lambda\|}{\|\lambda\|}
\widetilde{\mathcal{F}}\hspace{-.25mm}g(\lambda)\,.
\end{equation}
Now, by substituting \eqref{SolutionFourier}
and \eqref{SolutionFourierBis} in \eqref{PotentialFourier1},
we get
\begin{equation}\label{PotentialFourier2}
\begin{aligned}
\mathcal{P}[u](t)
&={\textstyle\frac{c_0}2}{\displaystyle\int_{\mathfrak{a}}}\;
\|\lambda\|^2\,(\cos t\|\lambda\|)^2\,
\mathcal{F}\hspace{-.25mm}f(\lambda)\,
\widetilde{\mathcal{F}}\!f(\lambda)\,
\nu(\lambda)\,d\lambda\\
&+{\textstyle\frac{c_0}2}{\displaystyle\int_{\mathfrak{a}}}\;
(\sin t\|\lambda\|)^2\,
\mathcal{F}\hspace{-.25mm}g(\lambda)\,
\widetilde{\mathcal{F}}\hspace{-.25mm}g(\lambda)\,
\nu(\lambda)\,d\lambda\\
&+{\textstyle\frac{c_0}4}{\displaystyle\int_{\mathfrak{a}}}\;
\|\lambda\|\,(\sin 2t\|\lambda\|)\,
\bigl\{\,\mathcal{F}\hspace{-.25mm}f(\lambda)\,
\widetilde{\mathcal{F}}\hspace{-.25mm}g(\lambda)
+\mathcal{F}\hspace{-.25mm}g(\lambda)\,
\widetilde{\mathcal{F}}\!f(\lambda)\,\bigr\}\,
\nu(\lambda)\,d\lambda\,.\\
\end{aligned}
\end{equation}
Similarly to $\mathcal P[u]$, we can rewrite the kinetic energy as
\begin{equation*}
\mathcal{K}[u](t)
=\,{\textstyle\frac{c_0}2}{\displaystyle\int_{\mathfrak{a}}}\;
\partial_{t}\mathcal{F}\hspace{-.25mm}u(t,\lambda)\,
\partial_{t}\widetilde{\mathcal{F}}\hspace{-.25mm}u(t,\lambda)\,
\nu(\lambda)\,d\lambda\,.
\end{equation*}
Using the following facts
\begin{equation*}
\begin{cases}
\;\partial_{t}\mathcal{F}\hspace{-.25mm}u(t,\lambda)
=-\,\|\lambda\|\,(\sin t\|\lambda\|)\,
\mathcal{F}\hspace{-.25mm}f(\lambda)
+(\cos t\|\lambda\|)\,
\mathcal{F}\hspace{-.25mm}g(\lambda)\,,\\
\;\partial_{t}\widetilde{\mathcal{F}}\hspace{-.25mm}u(t,\lambda)
=-\,\|\lambda\|\,(\sin t\|\lambda\|)\,
\widetilde{\mathcal{F}}\!f(\lambda)
+(\cos t\|\lambda\|)\,
\widetilde{\mathcal{F}}\hspace{-.25mm}g(\lambda)\,,\\
\end{cases}
\end{equation*}
we deduce that
\begin{equation}
\label{KineticFourier1}
\begin{aligned}
\mathcal{K}[u](t)
&={\textstyle\frac{c_0}2}
{\displaystyle\int_{\mathfrak{a}}}\;
\|\lambda\|^2\,(\sin t\|\lambda\|)^2\,
\mathcal{F}\hspace{-.25mm}f(\lambda)\,
\widetilde{\mathcal{F}}\!f(\lambda)\,
\nu(\lambda)\,d\lambda\\
&+{\textstyle\frac{c_0}2}{\displaystyle\int_{\mathfrak{a}}}\;
(\cos t\|\lambda\|)^2\,
\mathcal{F}\hspace{-.25mm}g(\lambda)\,
\widetilde{\mathcal{F}}\hspace{-.25mm}g(\lambda)\,
\nu(\lambda)\,d\lambda\\
&-{\textstyle\frac{c_0}4}{\displaystyle\int_{\mathfrak{a}}}\;
\|\lambda\|\,(\sin 2t\|\lambda\|)\,
\bigl\{\,\mathcal{F}\hspace{-.25mm}f(\lambda)\,
\widetilde{\mathcal{F}}\hspace{-.25mm}g(\lambda)
+\mathcal{F}\hspace{-.25mm}g(\lambda)\,
\widetilde{\mathcal{F}}\!f(\lambda)\,\bigr\}\,
\nu(\lambda)\,d\lambda\,.\\
\end{aligned}
\end{equation}
By suming up \eqref{PotentialFourier2} and \eqref{KineticFourier1},
we obtain   the conservation of the total energy\,:
\begin{equation*}
\mathcal{E}[u](t)
=\,{\textstyle\frac{c_0}2}
{\displaystyle\int_{\mathfrak{a}}}\;
\bigl\{\|\lambda\|^2\,
\mathcal{F}\!f(\lambda)\,
\widetilde{\mathcal{F}}\!f(\lambda)
+\mathcal{F}\hspace{-.4mm}g(\lambda)\,
\widetilde{\mathcal{F}}\hspace{-.25mm}g(\lambda)
\bigr\}\,\nu(\lambda)\,d\lambda
=\,\mathcal{E}[u](0)\,.
\end{equation*}
That is $\mathcal{E}[u](t)$ is independent of $t.$
\section{Equipartition of energy}
\noindent
This section is devoted to the proof
of \eqref{StrictEquipartition} and \eqref{AsymptoticEquipartition}.
Using the classical trigonometric identities for double angles, we can
rewrite the identities  \eqref{PotentialFourier2} and \eqref{KineticFourier1}
respectively as
\begin{equation*}
\begin{aligned}
\mathcal{P}[u](t)
&={\textstyle\frac{c_0}4}
{\displaystyle\int_{\mathfrak{a}}}\,
\bigl\{\|\lambda\|^2\,
\mathcal{F}\!f(\lambda)\,
\widetilde{\mathcal{F}}\!f(\lambda)
+\mathcal{F}\hspace{-.4mm}g(\lambda)\,
\widetilde{\mathcal{F}}\hspace{-.25mm}g(\lambda)
\bigr\}\,\nu(\lambda)\,d\lambda\\
&+{\textstyle\frac{c_0}4}
{\displaystyle\int_{\mathfrak{a}}}\;
(\cos2t\|\lambda\|)\,
\bigl\{\|\lambda\|^2
\mathcal{F}\!f(\lambda)\,
\widetilde{\mathcal{F}}\!f(\lambda)
-\mathcal{F}\hspace{-.4mm}g(\lambda)\,
\widetilde{\mathcal{F}}\hspace{-.25mm}g(\lambda)
\bigr\}\,\nu(\lambda)\,d\lambda\\
&+{\textstyle\frac{c_0}4}
{\displaystyle\int_{\mathfrak{a}}}\;
\|\lambda\|\,(\sin 2t\|\lambda\|)\,
\bigl\{\,\mathcal{F}\hspace{-.25mm}f(\lambda)\,
\widetilde{\mathcal{F}}\hspace{-.25mm}g(\lambda)
+\mathcal{F}\hspace{-.25mm}g(\lambda)\,
\widetilde{\mathcal{F}}\!f(\lambda)\,\bigr\}\,
\nu(\lambda)\,d\lambda\\
\end{aligned}
\end{equation*}
and
\begin{equation*}
\begin{aligned}
\mathcal{K}[u](t)
&={\textstyle\frac{c_0}4}
{\displaystyle\int_{\mathfrak{a}}}\,
\bigl\{\|\lambda\|^2\,
\mathcal{F}\!f(\lambda)\,
\widetilde{\mathcal{F}}\!f(\lambda)
+\mathcal{F}\hspace{-.4mm}g(\lambda)\,
\widetilde{\mathcal{F}}\hspace{-.25mm}g(\lambda)
\bigr\}\,\nu(\lambda)\,d\lambda\\
&-{\textstyle\frac{c_0}4}
{\displaystyle\int_{\mathfrak{a}}}\;
(\cos2t\|\lambda\|)\,
\bigl\{\|\lambda\|^2
\mathcal{F}\!f(\lambda)\,
\widetilde{\mathcal{F}}\!f(\lambda)
-\mathcal{F}\hspace{-.4mm}g(\lambda)\,
\widetilde{\mathcal{F}}\hspace{-.25mm}g(\lambda)
\bigr\}\,\nu(\lambda)\,d\lambda\\
&-{\textstyle\frac{c_0}4}
{\displaystyle\int_{\mathfrak{a}}}\;
\|\lambda\|\,(\sin 2t\|\lambda\|)\,
\bigl\{\,\mathcal{F}\hspace{-.25mm}f(\lambda)\,
\widetilde{\mathcal{F}}\hspace{-.25mm}g(\lambda)
+\mathcal{F}\hspace{-.25mm}g(\lambda)\,
\widetilde{\mathcal{F}}\!f(\lambda)\,\bigr\}\,
\nu(\lambda)\,d\lambda\,.\\
\end{aligned}
\end{equation*}
Hence
\begin{equation}
\label{Difference1}
\begin{aligned}
&\mathcal{P}[u](t)\hspace{-.3mm}-\hspace{-.2mm}\mathcal{K}[u](t)\,=\\
&={\textstyle\frac{c_0}2}
{\displaystyle\int_{\mathfrak{a}}}\,
(\cos2t\|\lambda\|)\,
\bigl\{\|\lambda\|^2
\mathcal{F}\!f(\lambda)\,
\widetilde{\mathcal{F}}\!f(\lambda)
-\mathcal{F}\hspace{-.4mm}g(\lambda)\,
\widetilde{\mathcal{F}}\hspace{-.25mm}g(\lambda)
\bigr\}\,\nu(\lambda)\,d\lambda\\
&+\,{\textstyle\frac{c_0}2}
{\displaystyle\int_{\mathfrak{a}}}\;
\|\lambda\|\,(\sin 2t\|\lambda\|)\,
\bigl\{\,\mathcal{F}\hspace{-.25mm}f(\lambda)\,
\widetilde{\mathcal{F}}\hspace{-.25mm}g(\lambda)
+\mathcal{F}\hspace{-.25mm}g(\lambda)\,
\widetilde{\mathcal{F}}\!f(\lambda)\,\bigr\}\,
\nu(\lambda)\,d\lambda\,.\\
\end{aligned}
\end{equation}
Introducing polar coordinates in $\mathfrak{a}$,
\eqref{Difference1} becomes
\begin{equation}
\label{Difference2}
\mathcal{P}[u](t)\hspace{-.3mm}-\hspace{-.2mm}\mathcal{K}[u](t)
={\textstyle\frac{c_0}2}
{\displaystyle\int_{\,0}^{+\infty}}\hspace{-1mm}
\bigl\{\cos(2tr)\,\Phi(r)+\sin(2tr)\,r\,\Psi(r)\bigr\}\,
r^{d-1}\,dr\,,
\end{equation}
where
\begin{equation*}
\begin{aligned}
\Phi(r)\,&=\int_{S(\mathfrak{a})}
\textstyle\bigl\{\,r^2\mathcal{F}\!f(r\sigma)\,
\widetilde{\mathcal{F}}\!f(r\sigma)
-\mathcal{F}\hspace{-.4mm}g(r\sigma)\,
\widetilde{\mathcal{F}}\hspace{-.4mm}g(r\sigma)\bigr\}\,
\nu(r\sigma)\,d\sigma\,,\\
\Psi(r)\,&=\int_{S(\mathfrak{a})}
\textstyle\bigl\{\mathcal{F}\!f(r\sigma)\,
\widetilde{\mathcal{F}}\hspace{-.4mm}g(r\sigma)
+\mathcal{F}\hspace{-.4mm}g(r\sigma)\,
\widetilde{\mathcal{F}}\!f(r\sigma)\bigr\}\,
\nu(r\sigma)\,d\sigma\,,\\
\end{aligned}
\end{equation*}
and $d\sigma$ denotes the surface measure
on the unit sphere $S(\mathfrak{a})$ in $\mathfrak{a}$.
Let \,$\gamma_0\!\in\!(0,+\infty]$ \,be the width of
the largest horizontal strip \,$|\Im z|\!<\!\gamma_0$
\,in which \,$z\mapsto\nu(z\sigma)$ \,is holomorphic
for all directions \,$\sigma\!\in\!S(\mathfrak{a})$.

\begin{lemma}
\label{EstimatePhiPsi}
{\rm (i)} $\Phi(z)$ and $\Psi(z)$ extend to even holomorphic functions
in the strip \,$|\Im z|\!<\!\gamma_0$.

{\rm (ii)} If $\gamma_0\!<\!+\infty$,
the following estimate holds in every substrip \,$|\Im z|\!\le\!\gamma$
\,with \,$\gamma\!<\!\gamma_0:$
\,For every $N\!>\!0$, there is a constant \,$C\!>\!0$
$($depending on $f,g\!\in\!\mathcal{C}_R^\infty(\mathfrak{a})$,
$N$ \!and $\gamma$$)$ such that
\begin{equation*}
|\Phi(z)|+|\Psi(z)|\le
C\,|z|^{|\mathcal{R}_0^+|}\,(1\!+\!|z|)^{-N}\,e^{\,2R\,|\Im z|}\,.
\end{equation*}

{\rm (iii)} If $\gamma_0=+\infty$,
the previous estimate holds uniformly in $\mathbb{C}$.
\end{lemma}

\begin{proof}
(i) follows from the definition of \,$\Phi$ and $\Psi$.
Let us turn to the estimates (ii) and (iii).
On one hand,
according to the Paley-Wiener Theorem
(Theorem \ref{PaleyWiener}), all transforms
\,$\mathcal{F}\!f(z\sigma)$,
$\widetilde{\mathcal{F}}\!f(z\sigma)$,
$\mathcal{F}\hspace{-.25mm}g(z\sigma)$,
$\widetilde{\mathcal{F}}\hspace{-.25mm}g(z\sigma)$
\,are \,$\text{O}\bigl(\{1\!+\!|z|\}^{-N}e^{\,R\,|\Im z|}\bigr)$.
On the other hand,
let us discuss the behavior of the Plancherel measure.
Consider first the case where all multiplicities are integers.
Without loss of generality, we may assume that $k_\alpha\!\in\!\mathbb{N}^*$
and $k_{2\alpha}\!\in\!\mathbb{N}$
\;for every indivisible root $\alpha$.
Then
\begin{equation*}
\begin{aligned}
\nu(\lambda)\,
&\textstyle
=\,\const\;\prod_{\alpha\in\mathcal{R}_0^+}
\langle\lambda,\check\alpha\rangle\,
\bigl\{\langle\lambda,\check\alpha\rangle\!
+\!i(k_{\alpha}\hspace{-.75mm}+\!2k_{2\alpha})\bigr\}\\
&\textstyle\times\;
\prod_{0<j<k_\alpha}
\bigl\{\langle\lambda,\check\alpha\rangle^2\hspace{-.75mm}+\!j^2\bigr\}
\;\prod_{0\le\widetilde{j}<k_{2\alpha}}
\bigl\{\langle\lambda,\check\alpha\rangle^2\hspace{-.75mm}
+\!(k_{\alpha}\hspace{-.75mm}+\!2\widetilde{j})^2\bigr\}
\end{aligned}
\end{equation*}
is a polynomial of degree
\,$2\,|k|=2\sum_{\alpha\in\mathcal{R}^+}\!k_\alpha$.
In general,
\begin{equation*}
\nu(\lambda)=\const\,\pi(\lambda)\,\widetilde\nu(\lambda)\,,
\end{equation*}
where
\begin{equation*}
\textstyle
\pi(\lambda)
=\,\prod_{\alpha\in\mathcal{R}_0^+}
\langle\lambda,\check\alpha\rangle
\end{equation*}
is a homogeneous polynomial of degree $|\mathcal{R}_0^+|$
and
\begin{equation*}
\textstyle
\widetilde\nu(\lambda)\,
=\hspace{-1mm}\prod\limits_{\alpha\in\mathcal{R}_0^+}\hspace{-1mm}
\frac
{\Gamma\bigl(i\langle\lambda,\check\alpha\rangle\,
+\,k_\alpha\bigr)\vphantom{\big|}}
{\Gamma\bigl(i\langle\lambda,\check\alpha\rangle\,+\,1\bigr)
\vphantom{\big|}}\,
\frac
{\Gamma\bigl(\frac{i \langle\lambda,\check\alpha\rangle\,
+\,k_{\alpha}}2+\,k_{2\alpha}\bigr)\vphantom{\big|}}
{\Gamma\bigl(\frac{i\langle\lambda,\check\alpha\rangle\,
+\,k_\alpha}2\bigr)\vphantom{\big|}}\,
\frac
{\Gamma\bigl(-\,i\langle\lambda,\check\alpha\rangle\,
+\,k_\alpha\bigr)\vphantom{\big|}}
{\Gamma\bigl(-\,i\langle\lambda,\check\alpha\rangle\,+\,1\bigr)
\vphantom{\big|}}\,
\frac
{\Gamma\bigl(\frac{-i \langle\lambda,\check\alpha\rangle\,
+\,k_{\alpha}}2+\,k_{2\alpha}+\,1\bigr)\vphantom{\big|}}
{\Gamma\bigl(\frac{-i\langle\lambda,\check\alpha\rangle\,
+\,k_\alpha}2\bigr)\vphantom{\big|}}
\end{equation*}
is an analytic function which never vanishes on $\mathfrak{a}$.
Notice that \,$z\mapsto\nu(z\sigma)\text{ or }\widetilde\nu(z\sigma)$
has poles for generic directions $\sigma\!\in\!S(\mathfrak{a})$
as soon as some multiplicities are not integers.
Using Stirling's formula
\begin{equation*}
\Gamma(\xi)\sim\sqrt{2\pi}\,\xi^{\xi-\frac12}\,e^{-\xi}
\quad\text{as \,}|\xi|\to+\infty
\text{ \,with \,}|\arg z|\!<\!\pi\!-\!\varepsilon\,,
\end{equation*}
we get the following estimate for the Plancherel density, in each
strip \,$|\Im z|\!<\!\gamma$ \,with \,$0\!<\!\gamma\!<\!\gamma_0$\,:
\begin{equation*}
|\nu(z\sigma)|\le\,C\,|z|^{|\mathcal{R}_0^+|}\,(1\!+\!|z|)^{2|k|-|\mathcal{R}_0^+|}\,.
\end{equation*}
The estimates (ii) and (iii) follow easily from these considerations.
\end{proof}

\begin{proposition}
\label{DifferenceEstimate1}
Assume that the dimension $d$ is odd
and that all multiplicities are integers.
Then there exists a constant \,$C\!>\!0$
$($depending on the initial data
\linebreak
$f,g\!\in\!\mathcal{C}_R^\infty(\mathfrak{a}))$
such that, for every \,$\gamma\!\ge\!0$ and \,$t\!\in\!\mathbb{R}$,
\begin{equation*}
|\mathcal{P}[u](t)\hspace{-.3mm}-\hspace{-.2mm}\mathcal{K}[u](t)|
\le C\,e^{\,2\gamma(R-|t|)}\,.
\end{equation*}
\end{proposition}

\begin{proof}
Evenness allows us to rewrite \eqref{Difference2} as follows\,:
\begin{equation*}
\mathcal{P}[u](t)\hspace{-.3mm}-\hspace{-.2mm}\mathcal{K}[u](t)
={\textstyle\frac{c_0}4}
{\displaystyle\int_{-\infty}^{+\infty}}\hspace{-1mm}e^{\,i2tr}
\bigl\{\Phi(r)\!-\!ir\Psi(r)\bigr\}\,r^{d-1}\,dr\,.
\end{equation*}
Let us shift the contour of integration
from \,$\mathbb{R}$ \,to \,$\mathbb{R}\hspace{-.75mm}\pm\!i\gamma$,
according to the sign of~$t$,
and estimate the resulting integral, using Lemma \ref{EstimatePhiPsi}.iii.
As a result, the difference of energy
$$\mathcal{P}[u](t)-\mathcal{K}[u](t)= \frac{c_0}4\,e^{-2\gamma|t|}
{\int_{-\infty}^{+\infty}}
e^{\,i2tr}\,\bigl\{\Phi(r\!\pm\!i\gamma)\!
-i(r\!\pm\!\gamma)\,\Psi(r\!\pm\!i\gamma)\bigr\}\,
(r\!\pm\!i\gamma)^{d-1} dr$$
is \,$\text{O}\bigl((1\!+\!\gamma)^{-N}e^{\,2\gamma(R-|t|)}\bigr)$.
\end{proof}

As an immediate consequence of the above statement and in view of the fact that
$\gamma_0=\infty$ when  $k$ is integer valued,  we deduce the strict equipartition of energy \eqref{StrictEquipartition}
for \,$|t|\!\ge\!R,$ by letting $\gamma\!\to\!\infty.$

Henceforth, we will drop the above assumption on $k.$  By resuming the proof of Proposition \ref{DifferenceEstimate1}
and using Lemma \ref{EstimatePhiPsi}.ii
instead of Lemma \ref{EstimatePhiPsi}.iii,
we obtain the following result.

\begin{proposition}
\label{EstimateDifference2}
Assume that the dimension $d$ is odd.
Then, for every \,$0\!<\!\gamma\!<\!\gamma_0$,
there is a constant \,$C\!>\!0$
$(\text{depending on the initial data
$f,g\!\in\!\mathcal{C}_R^\infty(\mathfrak{a})$})$
such that
\begin{equation*}
|\mathcal{P}[u](t)\hspace{-.3mm}-\hspace{-.2mm}\mathcal{K}[u](t)|
\le C\,e^{-2\gamma|t|}
\quad\forall\;t\!\in\!\mathbb{R}\,.
\end{equation*}
\end{proposition}

As a corollary,
we obtain the asymptotic equipartition of energy
\eqref{AsymptoticEquipartition} in the odd dimensional case,
with an exponential rate of decay.
In the even dimensional case,
the expression \eqref{Difference2} cannot be handled by complex analysis
and we proceed differently.

\begin{proposition}
\label{EstimateDifference2}
Assume that the dimension $d$ is even.
Then there is a constant \,$C\!>\!0$
$(\text{depending on the initial data
$f,g\!\in\!\mathcal{C}_R^\infty(\mathfrak{a})$})$
such that
\begin{equation*}
|\mathcal{P}[u](t)\hspace{-.3mm}-\hspace{-.2mm}\mathcal{K}[u](t)|
\le C\,(1\!+\!|t|)^{-d-|\mathcal{R}_0^+|}
\quad\forall\;t\!\in\!\mathbb{R}\,.
\end{equation*}
\end{proposition}

\begin{proof}
The problem lies in the decay at infinity.
According to lemma \ref{EstimatePhiPsi},
$\Phi(r)$ and $\Psi(r)$ are divisible by $r^{D}$,
where $D\!=\!|\mathcal{R}_0^+|$.
Let us integrate \eqref{Difference2}
\,$d\!+\!D$ times by parts.
This way
\begin{equation*}
\int_{\,0}^{+\infty}\hspace{-1mm}
\cos(2tr)\,r^{d-1}
\underbrace{\Bigl\{
\int_{S(\mathfrak{a})}
\mathcal{F}\hspace{-.4mm}g(r\sigma)\,
\widetilde{\mathcal{F}}\hspace{-.4mm}g(r\sigma)\,
\nu(r\sigma)\,d\sigma
\Bigr\}}_{\widetilde\Phi(r)}
dr
\end{equation*}
becomes
\begin{equation*}
\pm\;{\textstyle\frac{1\text{ \,or \,}0}{(2t)^{d+D}}}\,
{\textstyle\frac{(d\,+D)\,!}{(D+1)\,!}}\,
\bigl({\textstyle\frac\partial{\partial r}}\bigr)^{D+1}
\widetilde\Phi(r)\big|_{r=0}\,
\pm\int_{\,0}^{+\infty}\hspace{-1mm}
{\textstyle\frac{\cos(2tr)\text{ or }\sin(2tr)}{(2t)^{d+D}}}\,
\bigl({\textstyle\frac\partial{\partial r}}\bigr)^{d+D}
\bigl\{r^{d-1}\,\widetilde\Phi(r)\bigr\}\,dr
\end{equation*}
which is \,$\text{O}\bigl(|t|^{-d-D}\bigr)$.
Similarly
\begin{equation*}
\int_{\,0}^{+\infty}\hspace{-1mm}
\cos(2tr)\,r^{d+1}\,
\Bigl\{\int_{S(\mathfrak{a})}\!
\mathcal{F}\!f(r\sigma)\,\widetilde{\mathcal{F}}\!f(r\sigma)\,
\nu(r\sigma)\,d\sigma\Bigr\}\,dr
=\,\text{O}\bigl(|t|^{-d-D-2}\bigr)
\end{equation*}
and
\begin{equation*}
\int_{\,0}^{+\infty}\hspace{-1mm}
\cos(2tr)\,r^{d}\,\Psi(r)\,dr
=\,\text{O}\bigl(|t|^{-d-D-1}\bigr).
\end{equation*}
This concludes the proof of Proposition \ref{EstimateDifference2}.
\end{proof}

As a corollary,
we obtain the asymptotic equipartition of energy
\eqref{AsymptoticEquipartition} in the even dimensional case,
with a polynomial rate of decay.

\begin{remark}
Our result may not be optimal.
In the $W$\hspace{-1mm}-invariant case,
one obtains indeed the rate of decay
\,$\text{\rm O}\bigl(\{1\!+\!|t|\}^{-d-2|\mathcal{R}_0^+\!|}\bigr)$
as in \cite{BOS}.
\end{remark}

\references \nextref{S}{Ben Said, S.}{\em  Huygens' principle for
the wave equation associated with the trigonometric Dunkl-Cherednik
operators} {Math. Research Letters {\bf 13} (2006), no. 1, 43--58 }

\nextref{SO} {Ben Said, S. and {\O}rsted, B.}{\em The wave equation
for Dunkl operators} { Math. (N.S.) {\bf 16} (2005), no. 3--4,
351--391 }

\nextref{BO}{Branson, T. and   \'Olafsson, G.}{\em Equipartition of energy for waves in symmetric space}
{J. Funct. Anal. {\bf 97} (1991), no. 2, 403--416}

\nextref{BOP1} {Branson, T. \'Olafsson, G. Pasquale, A.}{\em The
Paley--Wiener theorem and the local Huygens' principle for compact
symmetric spaces}{ Math. (N.S.) {\bf 16} (2005), no. 3--4, 393--428
}

\nextref{BOP2} {Branson, T. and   \'Olafsson, G.} {\em The
Paley-Wiener theorem for the Jacobi transform and the local Huygens'
principle for root systems with even multiplicities}{ Math. (N.S.)
16 (2005), no. 3--4, 429--442 }

\nextref{BOS} {Branson, T. \'Olafsson, G.Schlichtkrull, H.} {\em
Huygens' principle in Riemannian symmetric spaces} {Math. Ann. 301
(1995), no. 3, 445--462}

\nextref{KY} {El Kamel, J. Yacoub, C.} {\em Huygens' principle and
equipartition of energy for the modified wave equation associated to
a generalized radial Laplacian} {Ann. Math Blaise Pascal 12 (2005),
no. 1, 147--160}

\nextref{H} { Helgason, S.} {\em Geometric analysis on symmetric
spaces} {Math. Surveys Monographs 39, Amer. Math. Soc. (1994)}

\nextref{O1} { Opdam, E. M. }{\em Harmonic analysis for certain
representations of graded Hecke algebras} {Acta. Math. 175 (1995),
75--121}

\nextref{O2} { Opdam, E. M. }{\em Lecture notes on Dunkl operators
for real and complex reflection groups}{ Math. Soc. Japan Mem. 8
(2000) }

\lastpage
\end{document}